\documentclass[12pt]{amsart}
\usepackage{amsmath,amssymb,eucal,amsthm}
\usepackage{graphicx}

\def\lee{ \hbox{$\le$\kern-.8em\lower.61ex\hbox{$-$}} }
\def\gee{ \hbox{$\ge$\kern-.8em\lower.61ex\hbox{$-$}} }
\def\Re{{\rm Re}}
\def\Im{{\rm Im}\ }
\newcommand{\dis}{\displaystyle}

\def\C{\mathbb{C}}
\def\R{\mathbb{R}}
\def\Zed{\mathbb{Z}}
\def\Q{\mathbb{Q}}

\def\e{\varepsilon}

\def\Z2{\zeta_2(s,s)}

\def\v4{\vspace{.4cm}}

\begin{document}
\begin{center}{\bf On the distribution of the zeros of 
the Euler double zeta-function $\zeta_2(s,s)$} 
\v4
\v4

Kohji Matsumoto,\ Tomokazu Onozuka and Isao Wakabayashi
\end{center}
\v4

{\bf Abstract}  We investigate the distribution of zeros of the function $\Z2$, namely the function defined by putting $s_1=s_2=s$ in the Euler double zeta-function $\zeta_2(s_1,s_2)$. Our theorems provide refinements of observations, obtained numerically in a previous paper by the first author and M. Sh{\=o}ji, with rigorous theoretical proofs. 
\v4

\noindent
\section{Introduction}\label{sec-1}

The series
$$\zeta_2(s_1,s_2)=\sum_{n_1=1}^{\infty}\sum_{n_2=1}^{\infty}\frac1{n_1^{s_1}(n_1+n_2)^{s_2}}$$
converges absolutely in the domain $\{(s_1,s_2)\in\C^2\ |\ \Re\ (s_1+s_2)>2,\Re\ s_2>1\}$.  It can be continued analytically to the whole $\C^2$ as a meromorphic function, and is called {\it the Euler double zeta-function}.  Our ultimate aim is to know the zero set of this function in two complex variables, but it is difficult.  So, in this paper, we investigate, putting $s_1=s_2\ (=s)$, the zeros of the function $\Z2$ in one variable. 

In the present paper,
we shall determine the zero free region of $\Z2$ in its absolutely convergent domain $\{s\in\C\ |\ \Re\ s>1\}$, and give a lower estimate of the number of zeros of $\Z2$ in any strip given in the region $\{\Re\ s>1\}\backslash$\{the zero free region\} or $\{1/2<\Re\ s<1\}$. Further we shall show that under the Riemann hypothesis the function $\Z2$ has no zero in the region $\{s=\sigma+it\ |\ 0\leq \sigma \leq \frac12-\sigma(t)\}$ where $\sigma(t)$ is a certain continuous positive function with $\lim_{|t|\to \infty}\sigma(t)=0$.  We shall show that for any $\sigma_1<0$, there is a positive number $t_0$ and $\Z2$ has no zero in $\{\sigma_1\leq \Re s <0,\ |\Im s|\geq t_0\}$. 
\v4

We put $D(\sigma_1,\sigma_2)=\{s=\sigma+it\ |\ \sigma_1<\sigma<\sigma_2,\ t\in\R\}$, and denote by $N(\sigma_1,\sigma_2,T)$ the number of zeros of $\Z2$ in the region $\{\sigma_1<\Re\ s<\sigma_2,\ 0\leq\Im s\leq T\}$.

In 2014, the first author and Sh\={o}ji \cite{5} calculated numerically zeros of $\Z2$ for the first time.  
They proved (\cite[Proposition 2.1]{5}):

(i) There is no zero of $\Z2$ when $\sigma=\Re\ s$ is sufficiently large.  

They observed: 

(ii) There are many zeros of $\Z2$ in $D(0,1)$.

(iii) It seems that there are more zeros in $D(0,1/2)$ than in $D(1/2,1)$.

(iv) At least in the range of computation, there is no zero lying on the line $\{\sigma=1/2\}$.

(v) There are some zeros in the region $\{\sigma>1\}$.

(vi) There are some zeros in the region $\{\sigma<0\}$, but it seems that there is no zero outside the real axis when $|\sigma|$ is sufficiently large, and the number of zeros becomes few and few when $|t|$ becomes large.
\v4

In a subsequent article \cite{6},
the first author and Sh\={o}ji used the above computation in order to investigate the zero set of the Euler double zeta-function $\zeta_2(s_1,s_2)$ as a function in two variables.

Real zeros of $\zeta_2(s,s)$ were studied by Kamano \cite{Kamano} and also by \cite{5}, 
and then, the first author, Matsusaka and Tanackov \cite{4} investigated more closely the behavior of real zeros of $\Z2$ and more generally of the $n$-fold multiple zeta-functions with identical arguments $\zeta_n(s,\ldots,s)$.
\v4
 
In the well-known harmonic product formula
$$\zeta(s_1)\zeta(s_2)=\zeta_2(s_1,s_2)+\zeta_2(s_2,s_1)+\zeta(s_1+s_2),$$
we put $s_1=s_2\ (=s)$, then we obtain
$$\Z2=\frac12(\zeta(s)^2-\zeta(2s)). \eqno{(1.1)}$$
This is the starting point of the investigation in the present paper.
From (1.1) we find that, for $\Re\ s>1$,
$\Z2=0$ if and only if  
$$\frac{\zeta(s)^2}{\zeta(2s)}=\frac{\dis\prod_p\left(1-\frac1{p^{2s}}\right)}{\dis\prod_p\left(1-\frac1{p^s}\right)^2}=\prod_p\frac{1+\frac1{p^s}}{1-\frac1{p^s}}=1,\eqno{(1.2)}$$
where $p$ runs through all primes.  From Lemma 1 in Subsection 2.1, we have, for $s=\sigma+it\ (\sigma>1)$,
$$\max_{t\in \R}\left|\arg \frac{1+\frac1{p^s}}{1-\frac1{p^s}}\right| = \arg \frac{1+\frac{i}{p^\sigma}}{1-\frac{i}{p^\sigma}}. \eqno{(1.3)}$$
Denote by $p_n$ the $n$-th prime.
We define the following two functions in $\sigma>1$:

\noindent
{\bf Definition 1} $\dis\theta_n=\arg \frac{1+\frac{i}{p_n^\sigma}}{1-\frac{i}{p_n^\sigma}}
\qquad (n\geq 1),$\hfill (1.4)\\ 
%

\noindent
{\bf Definition 2} $\dis g(\sigma)=\theta_1-\sum_{n\geq 2}\theta_n.$\hfill {(1.5)}
\v4

We shall see that, by Proposition 1 in Subsection 2.2, observing the graph of $g(\sigma)$\ (Fig. 1), $g(\sigma)$ has a unique zero $\sigma_b=1.8018642\cdots$.
\v4

Now we list up our main results.

On a zero free region of $\Z2$ in its absolutely convergent domain $\{s \in\C\ |\ \Re\ s>1\}$, we have
\v4

\noindent
{\bf Theorem 1} $\quad$ $\Z2$ has no zero in $\{s\in\C\ |\ \Re\ s>\sigma_b\}$.
\v4

Theorem 1 is a refinement of the above (i) of \cite{5}.
The value $\sigma_b$ is best-possible, because
there are many zeros on the immediate left of the vertical line $\{\Re\ s=\sigma_b\}$ by Theorem 2 below.

On the number of zeros in the absolutely convergent domain, we have the following
quantitative estimate for (v) of \cite{5}.
\v4

\noindent
{\bf Theorem 2} $\quad$ For any $1<\sigma_0<\sigma_b$ and any sufficiently small $\delta>0$, we have $\dis\liminf_{T\to\infty}\frac1T N(\sigma_0-\delta,\sigma_0+\delta,T)>0$.
\v4

On the number of zeros in the half of the critical strip $\{1/2<\Re\ s<1\}$, we have
\v4

\noindent
{\bf Theorem 3} $\quad$ For any $\frac12<\sigma_0<1$ and any sufficiently small $\delta>0$.
we have $\dis\liminf_{T\to\infty}\frac1T N(\sigma_0-\delta,\sigma_0+\delta,T)>0$.
\v4

Theorem 3 gives a quantitative estimate for (ii) of \cite{5}.
Note that Theorem 3 is not new, since Nakamura and Pa\'nkowski \cite[Theorem 3]{8} in 2016 obtained the same result. But we give a proof in this paper, since our proof is different from theirs.

On the zeros in the other half of the critical strip $\{0\leq\Re\ s<1/2\}$, we have
\v4

\noindent
{\bf Theorem 4} $\quad$ Under the Riemann hypothesis, there is a constant $C>0$ such that $\zeta_2(s,s)$ has no zero in $\dis\left\{s=\sigma+it\ |\ 0\leq\sigma\leq\frac12-\frac C{\log\log(|t|+4)}\right\}$.
\v4

Even though Theorem 4 depends on the Riemann hypothesis, Theorem 3 and Theorem 4 show that the phenomenon is completely different on the left-hand side from the situation on the right-hand side of the critical line $\{\sigma=1/2\}$.  For example, Theorems 3 and 4 give, in a sense, a counter answer to (iii) of \cite{5}. But it is difficult to observe the difference of this phenomenon from \cite[Figure 1]{5}.  Moreover, by Theorem 4, it remains the possibility of existence of an infinite number of zeros of $\Z2$ on the immediate left of the critical line $\{\sigma=1/2\}$. One may observe from \cite[Figure 1]{5} that the zeros of $\Z2$ gradually approach the critical line as the imaginary parts of zeros become greater.

Theorems 3 and 4 do not answer (iv) of \cite{5}.
  
On the zeros in the domain $\{\Re\ s<0\}$, we have, without the Riemann hypothesis,
\v4

\noindent
{\bf Theorem 5} $\quad$ For any $\sigma_1<0$, there exists $t_0>0$ such that $\Z2$ has no zero in $\{s=\sigma+it\ |\ \sigma_1\leq \sigma<0,\ |t|\geq t_0\}$. 
\v4

Theorem 5 gives an answer to (vi) of \cite{5}.
 
\section{The zeros of $\Z2$ in the absolutely convergent domain}\label{sec-2}
\subsection{A preliminary}

First we prepare the following Lemma 1.    Note that, applying it to $r=\frac1{p^{\sigma}}$ with prime $p$ and $\sigma>1$, we easily obtain the equality (1.3) mentioned in the introduction.
\v4

\noindent
{\bf Lemma 1} For $0<r<1$ and $0\leq\varphi\leq \pi$, let $\dis \theta=\arg \frac{1+ri}{1-ri}, \theta'=\arg(1+re^{i\varphi}),\ \theta''=\arg\frac1{1-re^{i\varphi}}$ with $0<\theta<\pi/2,\ 0\leq \theta'<\pi/2,\  0\leq\theta''<\pi/2$.  Then $0\leq\theta'+\theta''\leq \theta$, where the second equality holds when $\varphi=\pi/2$. 
\v4

\noindent
{\it Proof}. We have $\tan\frac{\theta}2=r$, and $\dis\tan\theta=\frac{2\tan\frac{\theta}2}{1-\tan^2\frac{\theta}2}=\frac{2r}{1-r^2}$.
We also have $\dis\tan\theta'=\frac{r\sin\varphi}{1+r\cos\varphi}$ and $\dis\tan\theta''=\frac{r\sin\varphi}{1-r\cos\varphi}$.
Hence,
$$\tan(\theta'+\theta'')=\frac{\tan\theta'+\tan\theta''}{1-\tan\theta'\tan\theta''}=\frac{2r\sin\varphi}{1-r^2}.$$
Therefore $\tan(\theta'+\theta'')\leq\tan\theta$.  Noting that $\theta'+\theta''<\pi/2$ by continuity of the tangent function, we then have $\theta'+\theta''\leq\theta$, and the equality holds when $\varphi=\pi/2$.\qed 
\v4

\subsection{Image of the $\sigma$-line by the map $w=\log(\zeta(s)^2/\zeta(2s))$}

For a fixed $\sigma$, the vertical line $\{s\in\C\ |\ s=\sigma+it\ (t\in\R)\}$ will be called the $\sigma$-{\it line}.  We first investigate for $\sigma>1$ the image of the $\sigma$-line by the map $w=\log(\zeta(s)^2/\zeta(2s))$ using the similar method as in \cite[\S 11.1 - \S 11.6]{10}. 

Viewing (1.2) we define a function in $s$:
\v4

\noindent{\bf Definition 3} \ $\dis F(s)=\log\frac{\zeta(s)^2}{\zeta(2s)}=\sum_{n\geq1}\log\frac{1+\frac1{p_n^s}}{1-\frac1{p_n^s}}\ \ \ (\Re\ s>1)$,\hfill (2.1)\\
where $p_n$ denotes the $n$-th prime, and we take the branch of the logarithm to be real for real $s$. 
\v4

We denote by $U_{\sigma}$ the image of the $\sigma$-line by the map $w=F(s)$: 
$$U_{\sigma}=\{F(s)\ |\ s=\sigma+it\ (t\in\R)\}.$$
We put $r_n=\frac1{p_n^{\sigma}}\ (n\geq1)$.
Moreover, instead of the parameter $t$ we take an infinite number of independent parameters $\varphi_1,\varphi_2,\ldots$, and put
$$V_{\sigma}=\left\{\sum_{n\geq1}\log\frac{1+r_ne^{i\varphi_n}}{1-r_ne^{i\varphi_n}}\ (\varphi_n\in\R)\right\}.$$
We have evidently $U_{\sigma}\subset V_{\sigma}$.

We investigate the shape of $V_{\sigma}$.
For this purpose,
let $V_n\ (n=1,2,\ldots)$ be the image of the circle $\{|z|=r_n\}$ by the map $w=\log\frac{1+z}{1-z}$.  Then $V_n$ is a convex closed curve enclosing the point $w=0$ (see [1, p.87, foot note]), and symmetric with respect to both the coordinate axes. The curve $V_n$ passes the point $w=\rho_n:=\log\frac{1+r_n}{1-r_n}$.
  When $\varphi_n=\frac{\pi}2$, by (1.4) we have 
$$\log\frac{1+r_ne^{i\varphi_n}}{1-r_ne^{i\varphi_n}}=\log\frac{1+r_ne^{i\pi/2}}{1-r_ne^{i\pi/2}}=\log\frac{1+r_ni}{1-r_ni}=i\arg\frac{1+r_ni}{1-r_ni}=i\theta_n,$$ 
so $V_n$ passes also the point $w=i\theta_n$.
  
The set 
$$V_1+V_2=\{\log\frac{1+r_1e^{i\varphi_1}}{1-r_1e^{i\varphi_1}}+\log\frac{1+r_2e^{i\varphi_2}}{1-r_2e^{i\varphi_2}}\ |\ \varphi_1,\varphi_2\in\R \}$$
is a closed ring-shaped area enclosed by an outside convex closed curve and an inside convex closed curve (see \cite[p.88]{1}, \cite[\S9]{2}, \cite[\S11.5, \S11.6]{10}).  The set $V_1+V_2$ is symmetric with respect to both the coordinate axes.  The outside convex closed curve passes the points $w=\rho_1+\rho_2$ and $w=i(\theta_1+\theta_2)$, and the inside convex closed curve passes the points $w=\rho_1-\rho_2$ and $w=i(\theta_1-\theta_2)$. The point $w=0$ is inside of the inside convex closed curve.

If $\theta_1-\sum_{n\geq 2}\theta_n>0$ and $\rho_1-\sum_{n\geq 2}\rho_n>0$, then the set $V_{\sigma}=V_1+V_2+\cdots$ is a closed ring-shaped area, because $V_{\sigma}$ has the inside closed curve passing the points $i(\theta_1-\sum_{n\geq 2}\theta_n)$ and $\rho_1-\sum_{n\geq 2}\rho_n$. 
\v4

\noindent
{\bf Claim.} If $\theta_1-\sum_{n\geq 2}\theta_n\geq 0$, then $\rho_1-\sum_{n\geq 2}\rho_n>0$.
\v4

\noindent
{\it Proof.} 
We have 
\begin{align*}
\rho_n=\log\frac{1+r_n}{1-r_n}&=(r_n-\frac12r_n^2+\frac13r_n^3+\cdots)+(r_n+\frac12r_n^2+\frac13r_n^3+\cdots)\\
&=2(r_n+\frac13r_n^3+\frac15r_n^5+\cdots).
\end{align*}
By $\tan\frac{\theta_n}2=r_n$ we have $\theta_n=2\tan^{-1}r_n=2(r_n-\frac13r_n^3+\frac15r_n^5-\frac17r_n^7+\cdots)$.
Hence
$$\frac{\rho_n}{\theta_n}=\frac{1+\frac13r_n^2+\frac15r_n^4+\frac17r_n^6+\cdots}{1-\frac13r_n^2+\frac15r_n^4-\frac17r_n^6+\cdots}=1+\frac{2(\frac13r_n^2+\frac17r_n^6+\frac1{11}r_n^{10}+\cdots)}{1-\frac13r_n^2+\frac15r_n^4-\frac17r_n^6+\cdots}.$$
On the denominator here, we note that
$$
1-\frac{1}{3}r_n^2<1-\frac13r_n^2+\frac15r_n^4-\frac17r_n^6+\cdots<1.
$$
Hence
$$\frac{\rho_1}{\theta_1}=1+\frac{2(\frac13r_1^2+\frac17r_1^6+\frac1{11}r_n^{10}+\cdots)}{1-\frac13r_1^2+\frac15r_1^4-\cdots}>1+\frac23r_1^2,$$
and for $n\geq 2$ we have
$$\frac{\rho_n}{\theta_n}=1+\frac{2(\frac13r_n^2+\frac17r_n^6+\frac1{11}r_n^{10}+\cdots)}{1-\frac13r_n^2+\frac15r_n^4-\cdots}
<1+\frac23r_1^2\left(\frac{r_n}{r_1}\right)^2\frac{1+\frac3{7}r_n^4+\frac{3}{11}r_n^8+\cdots}{1-\frac13r_n^2}.$$
From the latter, noting $\frac{r_n}{r_1}=\frac{2^{\sigma}}{p_n^{\sigma}}\leq\left(\frac23\right)^{\sigma}<\frac23$,
we see that
$$\frac{\rho_n}{\theta_n}<1+\frac23r_1^2\left(\frac23\right)^2\frac{1+\frac3{7}(r_n^4+r_n^8+\cdots)}{1-\frac13r_n^2}<1+\frac23r_1^2 \quad (n\geq 2),$$
because
\begin{align*}
\frac{1+\frac3{7}(r_n^4+r_n^8+\cdots)}{1-\frac13r_n^2}
=\frac{1+\frac{3}{7}\cdot\frac{r_n^4}{1-r_n^4}}
{1-\frac{1}{3}r_n^2}
=\frac{1-\frac{4}{7}r_n^4}{(1-r_n^4)(1-\frac{1}{3}r_n^2)}
<\frac{1}{(1-\frac{1}{3^4})(1-\frac{1}{3^3})}
\end{align*}
which is obviously less than $(3/2)^2$.
Hence we have 
$$\theta_1-\sum_{n\geq 2}\theta_n<\frac1{1+\frac23r_1^2}(\rho_1-\sum_{n\geq 2}\rho_n).$$ 
Therefore, if $\theta_1-\sum_{n\geq 2}\theta_n\geq 0$, then $\rho_1-\sum_{n\geq 2}\rho_n>0$.\qed
\v4

This claim means the following. When $\sigma$ varies from values with $\theta_1-\sum_{n\geq 2}\theta_n>0$ to the value with $\theta_1-\sum_{n\geq 2}\theta_n=0$, then the inside convex closed curve of $V_{\sigma}$ collapses to a segment not on the imaginary axis but on the real axis.

Hence, by this claim, the set $V_{\sigma}$ has two different types of shape (i) and (ii) as follows (see \cite[\S 11.6]{10}):

(i) If $\theta_1-\sum_{n\geq 2}\theta_n>0$, then $V_{\sigma}$ is a closed ring-shaped area enclosed by two convex closed curves, and the inside closed curve of $V_{\sigma}$ passes the point $i(\theta_1-\sum_{n\geq 2}\theta_n)$,\\

(ii) if $\theta_1-\sum_{n\geq 2}\theta_n\leq 0$, then $V_{\sigma}$ is a closed region enclosed by a single convex closed curve. 
\v4 

In order to determine the border value of two types (i) and (ii), we investigate the graph of $g(\sigma)\ (\sigma>1)$ (see Definition 2).
\v4

\noindent
{\bf Proposition 1}\ \ There exists a unique $\sigma_m>1$ such that $g(\sigma)$ is increasing in $1<\sigma<\sigma_m$ and decreasing in $\sigma_m<\sigma$.  Numerically $\sigma_m=2.6362977\cdots$ and $g(\sigma_m)=0.158267316\cdots>0$.  We have $\lim_{\sigma\to 1+0}g(\sigma)=-\infty$ and $\lim_{\sigma\to \infty}g(\sigma)=0$.  Further there exists a unique $\sigma_b>1$ with $g(\sigma_b)=0$, and $g(\sigma)<0$ for $1<\sigma<\sigma_b$, and $g(\sigma)>0$ for $\sigma_b<\sigma$. Numerically $\sigma_b=1.8018642\cdots$.
\v4

\noindent
{\it Proof}. Recall that by (1.4) and (1.5), we have for $\sigma>1$ 
$$g(\sigma)=\theta_1-\sum_{n\geq 2}\theta_n,\ \ \ \ \theta_n=\arg \frac{1+\frac{i}{p_n^\sigma}}{1-\frac{i}{p_n^\sigma}}.$$
Since $\tan \frac{\theta_n}2 =1/p_n^{\sigma}$, we have $\theta_n=2\tan^{-1}\frac1{p_n^{\sigma}}$. Differentiating this by $\sigma$ we have 
$$\frac{d\theta_n}{d\sigma}=-2\frac{\log p_n}{p_n^{\sigma}+\frac1{p_n^{\sigma}}}.$$

Hence
$$g'(\sigma)=\frac{2\log 2}{{2^{\sigma}}+\frac1{2^{\sigma}}}    \left(\sum_{p\geq3}\frac{\log p}{\log2}\cdot\frac{2^{\sigma}+1/2^{\sigma}}{p^{\sigma}+1/p^{\sigma}}-1\right).$$
For every prime number $p\geq 3$, putting
$$h_p(\sigma)=\frac{2^{\sigma}+1/2^{\sigma}}{p^{\sigma}+1/p^{\sigma}},$$
we have
$$h_p'(\sigma)=\frac{(\log 2)(2^{\sigma}-1/2^{\sigma})(p^{\sigma}+1/p^{\sigma})-(\log p)(2^{\sigma}+1/2^{\sigma})(p^{\sigma}-1/p^{\sigma})}{(p^{\sigma}+1/p^{\sigma})^2}.$$
Since we easily see
$$\frac{2^{\sigma}+1/2^{\sigma}}{2^{\sigma}-1/2^{\sigma}}>\frac{p^{\sigma}+1/p^{\sigma}}{p^{\sigma}-1/p^{\sigma}},$$
we have $h_p'(\sigma)<0$. Hence $h_p(\sigma)$ is a decreasing function of $\sigma$. So
$$h(\sigma):=\sum_{p\geq3}\frac{\log p}{\log2}\cdot\frac{2^{\sigma}+1/2^{\sigma}}{p^{\sigma}+1/p^{\sigma}}$$
is a decreasing function of $\sigma$, and we have $\lim_{\sigma\to 1+0}h(\sigma)=\infty,\ \lim_{\sigma\to \infty}h(\sigma)=0$. Hence there exists a unique $\sigma_m>1$ with $h(\sigma_m)=1$. So we have $g'(\sigma_m)=0$, and $g'(\sigma)>0$ for $1<\sigma<\sigma_m$, and $g'(\sigma)<0$ for $\sigma_m<\sigma$.  By numerical computation we have $\sigma_m=2.6362977\cdots,\ g(\sigma_m)=0.158267316\cdots>0$. Further we have $\lim_{\sigma\to 1+0}g(\sigma)=-\infty,\ \lim_{\sigma\to \infty}g(\sigma)=0$.  Therefore, there exists a unique $\sigma_b>1$ with $g(\sigma_b)=0$. And we have $g(\sigma)<0$ for $1<\sigma<\sigma_b$, and $g(\sigma)>0$ for $\sigma_b<\sigma$.  By numerical computation we have $\sigma_b=1.8018642\cdots$. (We used Mathematica for the computation.)
\v4

Fig. 1  Graph of $g(\sigma)=\theta_1-\sum_{n=2}^{\infty}\theta_n$.

\vspace{-5.5cm}
\hspace{2.5cm}
\includegraphics[width=12cm,bb=0 0 640 480]{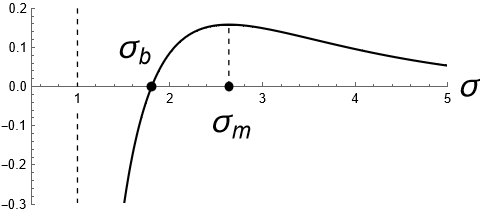}

\hfill$\Box$
\vspace{.5cm}

\noindent
{\bf Proposition 2}\ \ (i) If $\sigma>\sigma_b$, then $V_{\sigma}$ is a closed ring-shaped area enclosed by two convex closed curves.

(ii) If $1<\sigma<\sigma_b$, then $V_{\sigma}$ is a closed region enclosed by a convex closed curve.
\v4

\noindent
{\it Proof}.  The value $g(\sigma)=\theta_1-\sum_{n\geq 2}\theta_n>0$ if and only if $V_{\sigma}$ is a closed ring-shaped area enclosed by two convex closed curves, and the smaller closed curve of the boundary of $V_{\sigma}$ passes the point $i(\theta_1-\sum_{n\geq 2}\theta_n)$.  Hence if $\sigma>\sigma_b$, then $g(\sigma)>0$, and $V_{\sigma}$ is a closed ring-shaped area, and if $1<\sigma<\sigma_b$, then $g(\sigma)<0$, and $V_{\sigma}$ is not a closed ring-shaped area, i.e. $V_{\sigma}$ is a closed region enclosed by a single convex closed curve. \qed
\v4

The following lemma is necessary in the proof of Theorem 1 to determine the branch of logarithmic function.
\v4

\noindent
{\bf Lemma 2}\ \ For $\sigma>\sigma_b$, the ring-shaped area $V_{\sigma}$ does not contain the points $2\pi ni\ (n\in \Zed\backslash\{0\})$.
\v4

\noindent
{\it Proof}.  The outer convex closed  curve of $V_{\sigma}$ passes the point $i\sum_{n\geq 1}\theta_n$. We have
$$\frac{\theta_1}2=\arg\left(1+\frac i{2^{\sigma}}\right)<\arg\left(1+\frac i{2^{\sigma_b}}\right)<\frac1{2^{\sigma_b}},$$
so $\dis\theta_1<\frac1{2^{\sigma_b-1}}<1$. For $\sigma>\sigma_b$, we have $g(\sigma)=\theta_1-\sum_{n\geq 2}\theta_n>0$, so $\sum_{n\geq 1}\theta_n=\theta_1+\sum_{n\geq 2}\theta_n<\theta_1+\theta_1<2<\pi$.  Hence the points $2\pi ni\ (n\in \Zed\backslash\{0\})$ are not contained in $V_{\sigma}$.  \qed
\v4

\subsection{Proof of Theorem 1}

We fix $\sigma>\sigma_b$. By (1.2), the zeros of $\Z2$ coincide with the points satisfying $\frac{\zeta(s)^2}{\zeta(2s)}=1$, hence they coincide with the points satisfying $F(s)=\log\frac{\zeta(s)^2}{\zeta(2s)}=0$, since the possibility of $\log\frac{\zeta(s)^2}{\zeta(2s)}=2\pi ni$ with $n\in \Zed\backslash\{0\}$ is excluded by Lemma 2 and  $U_{\sigma}\subset V_{\sigma}$.  But by Proposition 2, $V_{\sigma}$ does not contain the point $w=0$, so we have 
$0\notin U_{\sigma}$.  Hence there is no zero of $\Z2$ on the $\sigma$-line.  Therefore, the region $\{s\in\C\ |\ \Re\ s>\sigma_b\}$ is a zero free region of $\Z2$.

\hfill(End of the proof of Theorem 1)
\v4

\subsection{Proof of Theorem 2}

Let us take any $\sigma_0$ with $1<\sigma_0<\sigma_b$ and $\delta>0$ with $1<\sigma_0-\delta,\  \sigma_0+\delta<\sigma_b$.  Since $1<\sigma_0<\sigma_b$, the set $V_{\sigma_0}$ is, by Proposition 2, a closed region enclosed by a convex closed curve, and $V_{\sigma_0}\ni 0$. Hence there exists a sequence of real numbers $\{\varphi_n^0\}_{n=1}^{\infty}$ such that $\dis\sum_{n\geq1}\log\frac{1+r_ne^{i\varphi_n^0}}{1-r_ne^{i\varphi_n^0}}=0\ \ (r_n=1/p_n^{\sigma_0})$. 

Based on the sequence $\{\varphi_n^0\}_{n=1}^{\infty}$, we shall find many  zeros of $\Z2$ in the strip $\sigma_0-\delta<\Re\ s\ (=\sigma)<\sigma_0+\delta$ by using the generalized Kronecker theorem below and Rouch{\'e}'s theorem.  

We define a function in the variable $s$:
\v4

\noindent
{\bf Definition 4}\ $\dis G(s)=\log\prod_{n\geq1}\frac{1+\frac{e^{i\varphi_n^0}}{p_n^s}}{1-\frac{e^{i\varphi_n^0}}{p_n^s}}\ \ \ (\Re\ s>1)$.\hfill(2.2)\\
\v4

We have $G(\sigma_0)=0$. Let $C$ be the circle with center $\sigma_0$ and radius $\delta$. Taking a smaller positive number $\delta$ if necessary, we may suppose that $G(s)\neq 0$  on the circle $C$.

We put $g=\min_{s\in C}|G(s)|$, and we
take any positive number $\e$ such that 
$$\e<\min\left\{\frac12,\ g,\ \ \frac {g}{16\pi\sum_{n\geq1}\frac1{p_n^{\sigma_0-\delta}}}\right\}.$$

We take then a sufficiently large integer $N$ such that we have, for any $s$ in the strip $\sigma_0-\delta\leq\Re\ s\leq\sigma_0+\delta$,  
$$\left|\log\prod_{n>N}\frac{1+\frac{e^{i\varphi_n^0}}{p_n^s}}{1-\frac{e^{i\varphi_n^0}}{p_n^s}}\right|<\frac{\e}4 \ \ {\rm and}\ \ \left|\log\prod_{n>N}\frac{1+\frac1{p_n^s}}{1-\frac1{p_n^s}}\right|<\frac{\e}4.\eqno{(2.3)}$$
\v4

The generalized Kronecker theorem is as follows (\cite[\S 11.7]{10}):
\v4

\noindent
{\bf The generalized Kronecker theorem} $\quad$ Let $\alpha_1,\ldots,\alpha_N$ be real numbers which are linearly independent over $\Q$, $\xi_1,\ldots,\xi_N$ be real numbers and $0<\gamma<\frac12$.  Then
$$\lim_{T\to\infty}\frac1T\mu\left\{t\in [0,T]\ |\ \max_{1\leq n\leq N}\parallel t\alpha_n-\xi_n\parallel<\gamma\right\}=(2\gamma)^N,$$
where $\parallel\cdot\parallel$ denotes the minimal distance to $\Zed$, and $\mu\{\cdot\}$ denotes the one- dimensional Lebesgue measure.
\v4

We apply this theorem to $\alpha_n=\frac1{2\pi}\log p_n,\ \xi_n=-\varphi_n^0/2\pi\ (1\leq n\leq N)$ and $\gamma=\e$. Then we have
$$\lim_{T\to \infty}\frac1T\mu\{t\in[0,T]\ |\ \max_{1\leq n\leq N}\parallel (t\log p_n+\varphi_n^0)/2\pi\parallel<\e\}=(2\e)^N.$$
Hence there exists a positive number $T_0$ such that for any $T>T_0$ we have
$$\mu\{t\in[0,T]\ |\ \max_{1\leq n\leq N}\parallel (t\log p_n+\varphi_n^0)/2\pi\parallel<\e\}>\frac12(2\e)^NT.\eqno{(2.4)}$$

Next we show that we can find at least $N_1=[\frac1{4\delta}(2\e)^NT]$ values $0\leq t_1<t_2<\cdots<t_{N_1}<T$ such that for all $1\leq j\leq N_1$ 
$$\max_{1\leq n\leq N}\parallel (t_j\log p_n+\varphi_n^0)/2\pi\parallel<\e,\eqno{(2.5)}$$
and such that $t_j+2\delta\leq t_{j+1}\ (1\leq j\leq N_1-1)$  (see \cite[\S 11.8]{10}, 
\cite[p.100]{1}). 
A proof of this fact is given as follows. Put
$$A=\{t\in[0,T]\ |\ \max_{1\leq n\leq N}\parallel (t\log p_n+\varphi_n^0)/2\pi\parallel<\e\}.$$
Then $\mu(A)>\frac12(2\e)^NT$ by (2.4).  In order to explain the proof, we  suppose for convenience that the set $A$ is a closed set even though it is not a closed set. First we take the smallest element $t_1$ of $A$. (Actually $A$ is not closed, so we take $t_1\in A$ sufficiently near to $\liminf A$, and we add similar changes in the following also.)  Since
$$\mu(A\cap[t_1+2\delta,T])=\mu(A\backslash[t_1,t_1+2\delta])\geq\mu(A)-2\delta>\frac12(2\e)^NT-2\delta,$$
we can take the minimal element $t_2\in A\cap[t_1+2\delta,T]$ if $\frac12(2\e)^NT-2\delta>0$. 
We can take in this way $0\leq t_1<t_2<\cdots<t_{N_1}<T$ with $t_j+2\delta\leq t_{j+1}$, since 
$$\frac12(2\e)^NT-(N_1-1)2\delta=\frac12(2\e)^NT-\left(\left[\frac1{4\delta}(2\e)^NT\right]-1\right)2\delta$$
$$\geq\frac12(2\e)^NT-\left(\frac1{4\delta}(2\e)^NT-1\right)2\delta=2\delta>0.$$

Now for each $t_j\ (1\leq j\leq N_1)$, we compare the values at $s$ of two functions (see Definitions 3, 4)
$$G(s)=\log\prod_{n\geq1}\frac{1+\frac{e^{i\varphi_n^0}}{p_n^s}}{1-\frac{e^{i\varphi_n^0}}{p_n^s}}\ \ \ {\rm and}\ \ \ F(s+it_j)=\log\prod_{n\geq1}\frac{1+\frac1{p_n^{s+it_j}}}{1-\frac1{p_n^{s+it_j}}}.$$
We estimate from above the finite sum
$$\left|\log\prod_{n\leq N}\frac{1+\frac{e^{i\varphi_n^0}}{p_n^s}}{1-\frac{e^{i\varphi_n^0}}{p_n^s}}-\sum_{n\leq N}\log\frac{1+\frac1{p_n^{s+it_j}}}{1-\frac1{p_n^{s+it_j}}}\right|,$$
which is, by changing the combination a little,
$$=\left|\sum_{n\leq N}\log\frac{1+p_n^{-s}e^{i\varphi_n^0}}{1+p_n^{-s-it_j}}+\sum_{n\leq N}\log\frac{1-p_n^{-s-it_j}}{1-p_n^{-s}e^{i\varphi_n^0}}\right|.$$

Using the inequalities $|e^{i\theta}-1|\leq|\theta|$ for real $\theta$, and $|\log(1+z_1)-\log(1+z_2)|\leq 2|z_1-z_2|$ for complex $z_1,z_2$ with $|z_1|,|z_2|<1$, from (2.5) we have, for $\Re\ s=\sigma$, 
$$\left|\log\frac{1+p_n^{-s}e^{i\varphi_n^0}}{1+p_n^{-s-it_j}}\right|<4\pi p_n^{-\sigma}\e,$$
and
$$\left|\log\frac{1-p_n^{-s-it_j}}{1-p_n^{-s}e^{i\varphi_n^0}}\right|<4\pi p_n^{-\sigma}\e.$$
Hence, for $\sigma_0-\delta<\sigma<\sigma_0+\delta$ we have

$$\left|\log\prod_{n\leq N}\frac{1+\frac{e^{i\varphi_n^0}}{p_n^s}}{1-\frac{e^{i\varphi_n^0}}{p_n^s}}-\sum_{n\leq N}\log\frac{1+\frac1{p_n^{s+it_j}}}{1-\frac1{p_n^{s+it_j}}}\right|<8\pi\left(\sum_{n\leq N}\frac1{p_n^{\sigma_0-\delta}}\right)\e.$$
Then we have, by (2.3), for $s\in C$,
$$|G(s)-F(s+it_j)|\leq 8\pi\left(\sum_{n\leq N}\frac1{p_n^{\sigma_0-\delta}}\right)\e+\frac{\e}4+\frac{\e}4<\frac g2+\frac{g}4+\frac{g}4=g.$$
Now we apply Rouch\'e's theorem to
$$F(s+it_j)=G(s)+(F(s+it_j)-G(s))\ \ (|s-\sigma_0|\leq \delta).$$
Since $G(\sigma_0)=0$ and $\max_{s\in C}|F(s+it_j)-G(s)|<g=\min_{s\in C}|G(s)|$, the function $F(s+it_j)$ has at least one zero in $\{|s-\sigma_0|<\delta\}$, i.e. there exists $s_j^0=\sigma_j^0+it_j^0$ with $|s_j^0-\sigma_0|<\delta$ such that $F(\sigma_j^0+it_j^0+it_j)=0$.  Since $\sigma_0-\delta<\sigma_j^0<\sigma_0+\delta$, this means that the equation $F(s)=0$ has solutions $\sigma_j^0+i(t_j^0+t_j)\ (1\leq j\leq N_1)$ in the strip $\{\sigma_0-\delta<\Re\ s<\sigma_0+\delta\}$.  
Moreover, if $j\neq j'$, then $t_j$ and $t_{j'}$ are apart at least $2\delta$ from each other, and from $-\delta<t_j^0<\delta$ and $-\delta<t_{j'}^0<\delta$,  we have $t_j^0+t_j\neq t_{j'}^0+t_{j'}$.  Hence $\sigma_j^0+i(t_j^0+t_j)\neq \sigma_{j'}^0+i(t_{j'}^0+t_{j'})$, i.e. the points $\sigma_j^0+i(t_j^0+t_j)\ (1\leq j\leq N_1)$ are all different.  Also $t_j^0+t_j<\delta+T$.   Since the greatest $t_{N_1}^0+t_{N_1}$ only might be greater than $T$, and the smallest $t_1^0+t_1$ only might be smaller than 0, so we throw away these two, and we keep the points $\sigma_j^0+i(t_j^0+t_j)\ (2\leq j\leq N_1-1)$, then for these points, we have $0<t_j^0+t_j<T$. Therefore, the equation $F(s)=0$ has at least $N_1-2$ solutions in $\{\sigma_0-\delta<\Re\ s<\sigma_0+\delta,\ 0<\Im s<T\}$, and we have 
$$\dis\liminf_{T\to\infty}\frac1T N(\sigma_0-\delta,\sigma_0+\delta,T)\geq \liminf_{T\to\infty}\frac1T(N_1-2)$$
$$\ \ \ \ \ =\liminf_{T\to\infty}\frac1T\left(\left[\frac1{4\delta}(2\e)^NT\right]-2\right)=\frac1{4\delta}(2\e)^N>0.$$
\hfill(End of the proof of Theorem 2)
\vspace{8mm}

\section{The zeros of $\Z2$ in the right half of the critical strip}

\subsection{The universality of $\log\zeta(s)$}

The universality theorem for the Riemann zeta-function
was first proved by Voronin \cite{11}, and Reich \cite{9} gave its convenient formulation.  
The universality theorem also holds for $\log\zeta(s)$ (cf. Karatsuba-Voronin \cite[Chap.VII]{3}).  We use it for the proof of Theorem 3.
\v4

\noindent
{\bf The universality theorem for $\log \zeta(s)$}\ \ Let $K$ be a compact set in the strip $\{1/2<\Re\ s<1\}$ with connected complement.  Let a function $f(s)$ be continuous on  $K$ and holomorphic in the interior of $K$.  Then, for any positive number $\e$, we have
$$\liminf_{T\to\infty}\frac1T\mu\{t\in [0,T]\ |\ \max_{s\in K} |\log\zeta(s+it)-f(s)|<\e\}>0,$$
where the branch of $\log\zeta(s+it)$ is taken so that it is real for real $s>1,t=0$ and is extended continuously along the intervals $[2,2+i(\Im s+t)],[2+i(\Im s+t),\Re\ s+i(\Im s+t)]$. 
\v4

\noindent
Remark. The statement of \cite[Chap.VII, Theorem 1]{3} is slightly weaker than the above, but the above theorem is in fact proved there. Moreover, it is implicitly assumed that for $t$ in the above theorem, the Riemann zeta-function does not encounter its possible zeros on the above intervals in order to define $\log\zeta(s+it)$.  
\v4

\subsection{Proof of Theorem 3}

We take any $\sigma_0$ with $\frac12<\sigma_0<1$, and a positive $\delta$ with $\frac12<\sigma_0-\delta,\ \ \sigma_0+\delta<1$.
Let us define a map:
\v4

\noindent
{\bf Definition 5} $\dis z=F_1(s)=\frac12\log\zeta(2s)\;\; (1/2<\Re\ s)$,
where we take the branch of the logarithm to be real for real $s$. 
\v4

\noindent
{\bf Claim} \ The image by $F_1$ of the strip $\{\sigma_0-\delta<\Re\ s<\sigma_0+\delta\}$ is contained in the rectangle
$$R=\ \left\{z=x+iy\ |\ \frac12\log\frac{\zeta(4(\sigma_0-\delta))}{\zeta(2(\sigma_0-\delta))}\leq x\leq\frac12\log\zeta(2(\sigma_0-\delta))\right.,$$
$$\hspace{6.5cm}\left. |y|\leq\frac{\pi}4\sum_p\frac1{p^{2(\sigma_0-\delta)}}\right\}.$$
\v4

\noindent
{\it Proof}.  Put $s=\sigma+it$ and $z=x+iy=F_1(s)=\frac12\log\zeta(2s)$.  We have $x=\frac12\log|\zeta(2s)|$ and $ y=\frac12\arg\zeta(2s)$. We first bound $x$.  From $\sigma_0-\delta<\sigma<\sigma_0+\delta$ and
$$\zeta(2s)=\frac{1}{\dis{\prod_p\left(1-\frac{e^{-i2t\log p}}{p^{2\sigma}}\right)}},$$
we have
$$|\zeta(2s)|<\frac{1}{\dis{\prod_p\left(1-\frac{1}{p^{2(\sigma_0-\delta)}}\right)}}=\zeta(2(\sigma_0-\delta)),$$
so $x=\frac12\log|\zeta(2s)|<\frac12\log\zeta(2(\sigma_0-\delta))$,
while we have      
$$|\zeta(2s)|>\frac{1}{\dis{\prod_p\left(1+\frac{1}{p^{2(\sigma_0-\delta)}}\right)}}=\frac{\zeta(4(\sigma_0-\delta))}{\zeta(2(\sigma_0-\delta))},$$
so $\dis x=\frac12\log|\zeta(2s)|>\frac12\log\frac{\zeta(4(\sigma_0-\delta))}{\zeta(2(\sigma_0-\delta))}$.
Hence the inequality in the definition of $R$ for its real part is sufficient.

We next bound $y$.  In general, for $0<r<1$, putting $\theta$ to be the argument between the real positive axis and the tangent line in the first quadrant from the origin to the circle with radius $r$ and center 1, we have $|\arg(1+re^{i\varphi})|\leq \theta$, and viewing the graph of $\sin^{-1}r$ we have $\sin^{-1}r<\frac{\pi}2 r$, so we have $|\arg(1+re^{i\varphi})|\leq \theta=\sin^{-1}r<\frac{\pi}2 r$.  From
$$\arg\zeta(2s)=-\sum_p\arg\left(1-\frac{e^{-i2t\log p}}{p^{2\sigma}}\right),$$
and from the inequality $|\arg(1+re^{i\varphi})|<\frac{\pi}2 r$ for $0<r<1$ and $\varphi\in\R$, and from $\sigma_0-\delta<\sigma<\sigma_0+\delta$ we have
$$|y|=\frac12|\arg\zeta(2s)|<\frac12\cdot\frac{\pi}2\sum_p \frac1{p^{2\sigma}}<\frac{\pi}4\sum_p\frac1{p^{2(\sigma_0-\delta)}}.$$
Hence the inequality in the definition of $R$ for its imaginary part is sufficient. 

\hfill $\Box$
\v4

We take any open disk $B_1$ which contains $R$ in its interior (see Fig. 2).  Let $\delta_1$ be the distance between $\partial B_1$(the boundary of $B_1$) and $R$.

Let $B=\{|s-\sigma_0|<\delta\}$, and define $z=\varphi(s)$ to be an affine mapping which maps $B$ just onto $B_1$. 

We now apply the universality theorem for $\log\zeta(s)$ to the holomorphic function $\varphi(s)$ on $\overline{B}$ (the closure of $B$).  Then, for any positive number $\e<\delta_1$, there is a positive constant $c(\e)$ such that we have
$$\liminf_{T\to\infty}\frac1T\mu\{t\in [0,T]\ |\ \max_{s\in \overline{B}} |\log\zeta(s+it)-\varphi(s)|<\e\}=c(\e)>0.$$

We set $\dis A=\{t\geq0\ |\ \max_{s\in \overline{B}}|\log\zeta(s+it)-\varphi(s)|<\e\}$.
There exists a positive number $T_0$ such that for any $T>T_0$ we have
$$\mu([0,T]\cap A)>\frac12 c(\e)T.\eqno{(3.1)}$$

We observe that inequality (3.1) obtained by the universality theorem for $\log\zeta(s)$ and inequality (2.4) in \S2.4 obtained by the generalized Kronecker theorem are similar.  So, similarly as in \S2.4, we can find at least $N_1=[\frac1{4\delta}c(\e)T]$ values $0\leq t_1<t_2<\cdots<t_{N_1}<T$ such that for all $1\leq j\leq N_1$
$$\max_{s\in \overline{B}} |\log\zeta(s+it_j)-\varphi(s)|<\e,\eqno{(3.2)}$$
and such that $t_j+2\delta\leq t_{j+1}\ (1\leq j\leq N_1-1)$.

Now we apply Rouch\'e's theorem to
$$\log\zeta(s+it_j)-\frac12\log\zeta(2(s+it_j))=(\varphi(s)-\frac12\log\zeta(2(s+it_j)))+(\log\zeta(s+it_j)-\varphi(s))$$
$$\ \ =:f(s)+h(s)\ \ \ (s\in \overline{B}).$$
First, if $s$ runs through $\partial\overline{B}$ counterclockwise once, then $\varphi(s)$ runs through $\partial\overline{B_1}$ counterclockwise once, while the point $F_1(s+it_j)=\frac12\log\zeta(2(s+it_j))$ remains in $R\subset B_1$, therefore, by the argument principle $f(s)$ has exactly one zero in $B$.  Moreover, we have $\min_{s\in\partial \overline{B}}|f(s)|\geq\delta_1$.  Second, by (3.2), we have $\max_{s\in\partial\overline{B}}|h(s)|<\e<\delta_1$.
Therefore, by Rouch\'e's theorem, the function $\log\zeta(s+it_j)-\frac12\log\zeta(2(s+it_j))$ has one zero $s_j^0=\sigma_j^0+it_j^0$ in $B$. So, $\zeta(\sigma_j^0+i(t_j^0+t_j))^2/\zeta(2(\sigma_j^0+i(t_j^0+t_j)))=1$. Since  $t_j\ (1\leq j\leq N_1)$ are each other apart from at least $2\delta$, the points $\sigma_j^0+i(t_j^0+t_j)\ (1\leq j\leq N_1)$ are all different.  
Hence, by (1.1) we have 
$$\zeta_2(\sigma_j^0+i(t_j^0+t_j),\sigma_j^0+i(t_j^0+t_j))=\frac12(\zeta(\sigma_j^0+i(t_j^0+t_j))^2-\zeta(2(\sigma_j^0+i(t_j^0+t_j))))=0,$$
and the equation $\Z2=0$ has at least $N_1-2$ zeros $\sigma_j^0+i(t_j^0+t_j)\ \  (2\leq j\leq N_1-1)$ in $\{\sigma_0-\delta<\Re\ s<\sigma_0+\delta,\ \ 0<\Im s<T\}$. Hence we obtain 
$$\dis\liminf_{T\to\infty}\frac1T N(\sigma_0-\delta,\sigma_0+\delta,T)\geq \liminf_{T\to\infty}\frac1T(N_1-2)$$
$$\ \ \ \ \ =\liminf_{T\to\infty}\frac1T\left([\frac1{4\delta} c(\e)T]-2\right)=\frac{c(\e)}{4\delta}>0.$$
\v4
\newpage
Fig. 2 
\begin{center}
\vspace{-0.5cm}\includegraphics[trim=4cm 18cm 6cm 3.5cm, clip]{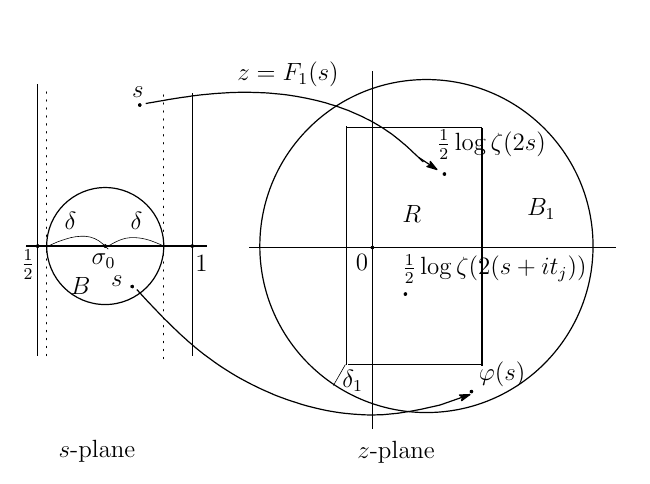}
\end{center}
\vspace{0.5cm}

\hfill(End of the proof of Theorem 3)
\vspace{8mm}

\section{The zeros of $\Z2$ in $\{0\leq\Re\ s<1/2\}$ under the Riemann hypothesis}
\label{sec-4}

\subsection{Order estimates}

In order to treat the zeros of $\Z2$ in $\{0\leq\Re\ s<1/2\}$, we need to use the order estimate of the Riemann zeta-function, i.e. the estimate of growth of $|\zeta(\sigma+it)|$ as a function in $t$ for fixed $\sigma$.  We prepare several known facts on the order estimation.

First we observe that
the process of the proof of the prime number theorem implies the following (cf. \cite[(3.11.8)]{10}):
\v4

\noindent
{\bf Assertion 1}  There is a constant $C>0$ such that 
$$\frac1{\zeta(s)}=O(\log (|t|+4))$$ 
uniformly in $\{\sigma\geq 1-\frac{C}{\log (|t|+4)},\ |t|\geq 1\}$.
\v4

\noindent
{\bf Assertion 2} (\cite[Theorem 13.18]{7})  Under the Riemann hypothesis, there is a constant $C_1>0$ such that $$|\zeta(s)|<\exp\left(\frac {C_1\log (|t|+4)}{\log\log(|t|+4)}\right)$$ for $\sigma\geq\frac12,\ |t|\geq 1.$
\v4

\noindent
{\bf  Assertion 3} (\cite[Theorem 13.23]{7})  Under the Riemann hypothesis, there is a constant $C_2>0$ such that  
$$\frac1{|\zeta(s)|}<\exp\left(\frac {C_2\log (|t|+4)}{\log\log(|t|+4)}\right)$$ for $|t|\geq 1,\ \sigma\geq\frac12+\frac1{\log\log(|t|+4)}.$
\v4

The Riemann hypothesis implies the Lindel{\"o}f hypothesis, so we have
\v4

\noindent
{\bf  Assertion 4}  Under the Riemann hypothesis, for any $\e>0$ we have $\zeta(s)=O(|t|^{1/2-\sigma+\e})$ uniformly in $\{0\leq\sigma\leq\frac12,\ |t|\geq 1\}$.
\v4

\noindent
{\bf  Assertion 5}\ (cf. \cite[Chap.~V]{10})  For any $\e>0$ and $\sigma_1<0$ we have $\zeta(s)=O(|t|^{1/2-\sigma+\e})$ uniformly in $\{\sigma_1\leq\sigma\leq 0,\ |t|\geq 1\}$.

\v4

\subsection{Proof of Theorem 4}

In this subsection we assume the truth of the Riemann hypothesis.  We use constants $C_1,\ C_2$ in Assertions 2 and 3.  Take a constant $C>\max\{C_1+C_2,1\}$.  The region in the statement of Theorem 4, restricting to $t\geq 1$, is
$$D_1=\left\{s=\sigma+it\ |\ t\geq 1,\ 0\leq \sigma\leq \frac12-\frac {C}{\log\log(t+4)}\right\}.$$

Under the Riemann hypothesis we first estimate $|\zeta(s)|$ for $s\in D_1$ from below. For this we use the functional equation
$$\zeta(s)=2^s\pi^{s-1}\Gamma(1-s)\sin\frac{\pi s}2\zeta(1-s).\eqno{(4.1)}$$
For $s\in D_1$ we have $\frac12+\frac {C}{\log\log(t+4)}\leq \Re\ (1-s)\leq 1$, hence from Assertion 3 we have in $D_1$ 
$$|\zeta(1-s)|>\exp\left(\frac {-C_2\log (t+4)}{\log\log(t+4)}\right).\eqno{(4.2)}$$
We use also the Stirling formula: For $\delta>0$,
$$\Gamma(s)=\sqrt{2\pi}s^{s-1/2}e^{-s}e^{O(|s|^{-1})}\ \ (|\arg s|\leq \pi-\delta,\ |s|\geq\delta),\eqno{(4.3)}$$
where the $O$-constant depends on $\delta$.  Then for $s\in D_1$ we have
$$|\Gamma(1-s)|=\sqrt{2\pi}|1-s|^{1/2-\sigma}e^{t\arg(1-s)}e^{\sigma-1}e^{O(|1-s|^{-1})},$$
and by $\arg(1-s)>-\pi/2$ for $s\in D_1$, and by taking the $O$-constant to be $C_{\delta}>0$, we have
$$|\Gamma(1-s)|>\sqrt{2\pi}|1-s|^{1/2-\sigma}e^{-t\pi/2}e^{\sigma-1}e^{-C_{\delta}}.\eqno{(4.4)}$$
Hence we have uniformly in $D_1$
$$|\Gamma(1-s)|\gg t^{1/2-\sigma}e^{-t\pi/2}.  \eqno{(4.5)}$$
For $s\in D_1$ we have also 
$$|\sin\frac{\pi s}2|\geq\frac12e^{t\pi/2}(1-e^{-t\pi}).\eqno{(4.6)}$$
Putting (4.2),(4.5) and (4.6) into (4.1) we then obtain uniformly in $D_1$
$$|\zeta(s)|\gg t^{1/2-\sigma}e^{\frac {-C_2\log (t+4)}{\log\log(t+4)}}.\eqno{(4.7)}$$

Next we estimate $|\zeta(2s)|$ for $s\in D_1$ from above, dividing into two cases i) and ii):

i) When $\Re\ s\geq \frac14$.  Since $\Re\ 2s\geq\frac12$, by Assertion 2 we have $$|\zeta(2s)|<e^{\frac {C_1\log (2t+4)}{\log\log(2t+4)}}.\eqno{(4.8)}$$

ii) When $0\leq\Re\ s<\frac14$.  Let $0<\e<\frac12$.  Since $0\leq\Re\ 2s<\frac12$, by Assertion 4 we have 
$$\zeta(2s)=O((2t)^{1/2-2\sigma+\e})\eqno{(4.9)}$$
uniformly in $\{0\leq\sigma<\frac14,\ t\geq 1\}$.

Thus we have the following:

i) When $\Re\ s\geq \frac14$. Since $\frac12-\sigma\geq \frac {C}{\log\log(t+4)}$ for $s\in D_1$, we have by (4.7) and (4.8) 
$$\frac{|\zeta(s)|^2}{|\zeta(2s)|}\gg t^{(1/2-\sigma)2}e^{\frac {-2C_2\log (t+4)}{\log\log(t+4)}}e^{\frac {-C_1\log (2t+4)}{\log\log(2t+4)}}$$
$$\geq e^{\frac {2C\log t}{\log\log(t+4)}}e^{\frac {-2C_2\log (t+4)}{\log\log(t+4)}}e^{\frac {-C_1\log (2t+4)}{\log\log(2t+4)}},$$
hence by $C>C_1+C_2$, we can find $t_0>0$ such that $\frac{|\zeta(s)|^2}{|\zeta(2s)|}>1$ for $t\geq t_0$. So by (1.1) there is no zero of $\Z2$ when  $t\geq t_0$ and $\frac14\leq \sigma\leq\frac12-\frac {C}{\log\log(t+4)}$.

ii) When $0\leq\Re\ s<\frac14$. By (4.7) and (4.9) we have
$$\frac{|\zeta(s)|^2}{|\zeta(2s)|}\gg t^{(1/2-\sigma)2}e^{\frac {-2C_2\log (t+4)}{\log\log(t+4)}}(2t)^{-(1/2-2\sigma+\e)}$$
$$\gg t^{1/2-\e}e^{\frac {-2C_2\log (t+4)}{\log\log(t+4)}},$$
hence we can find $t_0>0$ such that $\frac{|\zeta(s)|^2}{|\zeta(2s)|}>1$ for $t\geq t_0$, so by (1.1) there is no zero of $\Z2$ when  $t\geq t_0$ and $0\leq\sigma<\frac14$. 

At the beginning we took a constant $C$, and the results of two cases i) and ii) show that we find $t_0>0$ such that there is no zero of $\Z2$ in $D_1\cap \{t\geq t_0\}$, but by retaking $C$ sufficiently large we can eliminate this condition $t\geq t_0$, and we complete the proof of Theorem 4.\\
 \hfill (End of the proof of Theorem 4)
\vspace{8mm}

\section{The zeros of $\Z2$ in $\{\Re\ s<0\}$}

We prove Theorem 5.  For $\sigma_1<0$ let us consider $\Z2$ in the region $D_2=\{s=\sigma+it\ |\ \sigma_1\leq\sigma<0,\ t\geq 1\}$.

Using the functional equation (4.1), we first estimate $|\zeta(s)|$ for $s\in D_2$ from below. For $s\in D_2$ we have $1<\Re\ (1-s)\leq1-\sigma_1$, so by Assertion 1 we have uniformly in $D_2$
$$|\zeta(1-s)|\gg \frac1{\log(t+4)}.\eqno{(5.1)}$$
We have the same lower bound of $|\Gamma(1-s)|$ in $D_2$ as (4.5).  We have also the same lower bound of $|\sin\frac{\pi s}2|$ as (4.6).  Putting (5.1), (4.5) and (4.6) into (4.1) we then obtain uniformly in $D_2$
$$|\zeta(s)|\gg t^{1/2-\sigma}\frac1{\log(t+4)}.\eqno{(5.2)}$$

Next we estimate $|\zeta(2s)|$ for $s\in D_2$ from above.  Let $0<\e<\frac12$.  Since $2\sigma_1\leq\Re\ 2s<0$, by Assertion 5 we have uniformly in $D_2$
$$\zeta(2s)=O((2t)^{1/2-2\sigma+\e}).\eqno{(5.3)}$$
Hence by (5.2) and (5.3) we have uniformly in $D_2$
$$\left|\frac{\zeta(s)^2}{\zeta(2s)}\right|\gg t^{1/2-\e}\frac1{(\log(t+4))^2}.$$

Therefore, we find $t_0>0$ such that in $\{\sigma_1\leq \sigma<0,\ t\geq t_0\}$ we have $\left|\frac{\zeta(s)^2}{\zeta(2s)}\right|>1$. Hence by (1.1) in this region $\Z2$ has no zero.

\hfill (End of the proof of Theorem 5)
\vspace{8mm}

\noindent
{\bf Acknowledgment}  The authors express their gratitude to Dr.~Keita Nakai for his comment on the proofs in \S 3 and to Dr.~Yuichiro Toma for his comment on the proofs in \S 4.
The first author is supported by Japan Society for the Promotion of Science, Grant-in-Aid for
Scientific Research no. 22K03267.
\vspace{1.2cm}


Key words: the Euler double zeta-function, distribution of zeros
\v4

AMS subject classification: 11M32
\v4

Kohji Matsumoto, Graduate School of Mathematics, Nagoya University, Furocho, Chikusa-ku,
Nagoya 464-8602, Japan\\
and\\
Center for General Education, Aichi Institute of Technology, 1247 Yachigusa,
Yakusa-cho, Toyota 470-0392, Japan\\
E-mail address: kohjimat@math.nagoya-u.ac.jp

Tomokazu Onozuka, Faculty of Science and Technology, Oita University, 700 Dannoharu, Oita, 870-1192, Japan
\\
E-mail  address: math.onozuka@gmail.com

Isao Wakabayashi, Faculty of Science and Technology, Seikei University, Kichijoji Kitamachi, Musashino-shi, Tokyo, 180-8633, Japan\\
E-mail address: wakaba.isao@gmail.com

\end{document}